\documentclass[reqno,10pt]{amsart}
\usepackage{amssymb,amsmath}
\usepackage{mathrsfs}
\usepackage{latexsym}
\usepackage{amscd}
\usepackage{color}  
\usepackage{tikz}
\usepackage{caption}
\usepackage{subcaption}
\usetikzlibrary{arrows,calc,positioning,through,intersections,backgrounds,matrix,fit,shapes,decorations}
\usepackage[normalem]{ulem}  

\tikzset{
  vertex/.style={circle,color=black,draw=black,fill=black,inner sep=0cm,minimum size=.2cm},
  edge/.style={thick, black},
}

\let\cal=\mathcal

\newtheorem{theorem}{Theorem}[section]
\newtheorem{definition}[theorem]{Definition}

\newtheorem{prop}[theorem]{Proposition}

\newtheorem{cor}[theorem]{Corollary}

\newtheorem{claim}[theorem]{Claim}
\newtheorem{remark}[theorem]{Remark}

\def\epsilon{\varepsilon}

\def\cA{{\cal A}}
\def\cB{{\cal B}}
\def\cC{{\cal C}}

\def\cR{{\cal R}}

\def\rep{{\rm rep}}

\begin{document}

\title{On Odd Rainbow Cycles in Edge-Colored Graphs}

\author[A.~Czygrinow]{Andrzej Czygrinow}\thanks{The first author was partially supported
by Simons Foundation Grant $\#$521777.}    
\address{School of Mathematical and Statistical Sciences, 
Arizona State University, Tempe, AZ 85281, USA}
\email{aczygri@asu.edu}  

\author[T.~Molla]{Theodore Molla}\thanks{The second author was partially supported by 
NSF Grants DMS~1500121
and DMS~1800761.}  
\address{Department of Mathematics and Statistics, University of 
South Florida, Tampa, FL 33620, USA}
\email{molla@usf.edu}  

\author[B.~Nagle]{Brendan Nagle}\thanks{The third author was partially supported
by NSF Grant DMS~1700280.}  
\address{Department of Mathematics and Statistics, University of 
South Florida, Tampa, FL 33620, USA}  
\email{bnagle@usf.edu}  

\author[R.~Oursler]{Roy Oursler}
\address{School of Mathematical and Statistical Sciences, 
Arizona State University, Tempe, AZ 85281, USA}  
\email{Roy.Oursler@asu.edu}


\begin{abstract}
Let $G = (V, E)$ be an $n$-vertex edge-colored graph.  In~2013, H.~Li proved 
that 
if every vertex $v \in V$ is incident to at least $(n+1)/2$ distinctly colored edges, 
then $G$ admits a rainbow triangle.  
We prove that the same hypothesis ensures a rainbow $\ell$-cycle $C_{\ell}$ whenever 
$n \geq 432 \ell$.
This 
result is sharp for all odd integers $\ell \geq 3$, and extends earlier work of the
authors for when $\ell$ is even.    
\end{abstract}

\maketitle


\section{Introduction}

An {\it edge-colored graph} is a pair $(G, c)$, where $G = (V, E)$ is a graph and 
$c : E \to P$ is a function mapping edges to some palette of colors $P$.  
A subgraph $H \subseteq G$ is a {\it rainbow} subgraph if 
the edges of $H$ are distinctly
colored by $c$.   
Rainbow subgraph problems are a well-studied area of graph theory
(see, e.g., 
\cite{alon2017random}--\cite{LiWang}, and 
Section~\ref{sec:1.1} below).  
Here, 
we consider degree conditions on $(G, c)$  
ensuring the existence
of rainbow cycles $C_{\ell}$ of fixed length $\ell \geq 3$.   
To that end, 
a vertex $v \in V$  
in an edge-colored graph $(G, c)$ 
has {\it $c$-degree} $\deg_G^c(v)$ 
given by the number of distinct colors assigned by $c$ to the edges $\{v, w\} \in E$.
We   
set $\delta^c(G) = \min_{v \in V} \deg^c_G(v)$ for the minimum $c$-degree in $G$.  
The following result of H.~Li~\cite{li2013rainbow}
motivates our current work.  

\begin{theorem}[H.~Li~\cite{li2013rainbow}, 2013]
\label{thm:Li}  
Let $(G, c)$ be an $n$-vertex edge-colored graph.  
If $\delta^c(G) \geq (n+1)/2$, then 
$(G,c)$ admits a rainbow 3-cycle $C_3$.  
\end{theorem}  

\noindent   
A rainbow 
$K_{\lfloor n/2 \rfloor, \lceil n/2 \rceil}$ establishes that  
Theorem~\ref{thm:Li} is best possible.

We prove an analogue of 
Theorem~\ref{thm:Li} for $\ell$-cycles $C_{\ell}$ of fixed arbitrary length.

\begin{theorem}
\label{thm:odd}
For every integer $\ell \geq 3$, every edge-colored graph $(G, c)$ 
on $n \geq n_0(\ell)$
many vertices satisfying $\delta^c(G) \geq (n+1)/2$ admits a rainbow $\ell$-cycle 
$C_{\ell}$.  
\end{theorem}  

\noindent  A rainbow 
$K_{\lfloor n/2 \rfloor, \lceil n/2 \rceil}$ also establishes that  
Theorem~\ref{thm:odd} is best
possible for all odd integers $\ell$.

For even integers $\ell \geq 4$, 
the authors earlier proved in~\cite{CMNO} a stronger form of Theorem~\ref{thm:odd}.

\begin{theorem}[Czygrinow et al.~\cite{CMNO}]  
\label{thm:even}
For every even integer $\ell \geq 4$,  
every 
edge-colored graph
$(G, c)$ 
on $n \geq N_0(\ell)$ many vertices 
satisfying 
$\delta^c(G) \geq (n+5)/3$ 
admits a rainbow $\ell$-cycle $C_{\ell}$.  
\end{theorem}  

\noindent  It was shown in~\cite{CMNO} that 
Theorem~\ref{thm:even} is best
possible for every even $\ell \not\equiv 0$ (mod 3).

Theorem~\ref{thm:Li} holds non-vacuously when $n \geq 3$, and 
one may seek to quantify $n_0(\ell)$ and $N_0(\ell)$   
in Theorems~\ref{thm:odd} and~\ref{thm:even}.  
The proof of Theorem~\ref{thm:even}
depends on 
an application of the Szemer\'edi Regularity Lemma~\cite{Szem1, Szem2}, and therefore
gives very poor bounds on $N_0(\ell)$.  
Our proof of Theorem~\ref{thm:odd} is elementary, and easily provides
$n_0(\ell) = O(\ell^2)$.   
For the interested Reader, 
we provide a more detailed analysis in our final section 
which establishes
that $n_0(\ell)$ is linear in $\ell$.    
\begin{theorem}
\label{thm:ell}  
The function $n_0(\ell)$ in Theorem~\ref{thm:odd} satisfies $n_0(\ell) \le 432 \ell$.  
\end{theorem}  

The remainder of this paper is organized as follows.  
In Section~\ref{sec:1.1}, we discuss further results 
and context regarding 
rainbow cycle problems.  
In Section~\ref{sec:Li}, we sketch Li's proof~\cite{li2013rainbow}
of Theorem~\ref{thm:Li} and note the elements 
there 
which provide a basis for our approach here.  
In 
Section~\ref{sec:tools},   
we extend this proof to develop several tools useful for 
proving Theorems~\ref{thm:odd} and~\ref{thm:ell}.   
In Section~\ref{sec:odd}, we prove Theorem~\ref{thm:odd}, 
and in Section~\ref{sec:ell}, we prove Theorem~\ref{thm:ell}.  
In the entirety of this paper, we employ the following observations.  
\begin{remark}
\label{rmk:edgeminimal} 
\rm
We say that an edge-colored graph $(G, c)$ is {\it edge-minimal} when every 
$e \in E(G)$ satisfies $\delta^c(G - e) < \delta^c(G)$.  
Every edge-colored graph $(G, c)$  
admits an edge-minimal spanning 
subgraph $H \subseteq G$ satisfying $\delta^c(G) = \delta^c(H)$,       
so in Theorems~\ref{thm:Li}--\ref{thm:ell}  
it suffices to assume 
that $(G, c)$
is already edge-minimal.      
As such, $(G, c)$ admits no 
three commonly colored edges $\{u, v\}, \{v, w\}, \{w, x\} \in E(G)$, 
as removing $\{v, w\} \in E(G)$ violates edge-minimality.  
\hfill $\Box$  
\end{remark}  

\subsection{Rainbow cycles and anti-Ramsey theory}  
\label{sec:1.1}  
Li et al.~\cite{li2014rainbow}
extended Theorem~1.1 as follows:      
if the average   
$c$-degree $\alpha^c(G) = (1/n) \sum_{v \in V} \deg_G^c(v)$ 
satisfies 
$\alpha^c(G) \geq (n+1)/2$,  
then $(G, c)$ admits 
a rainbow $C_3$; 
if 
the number of colors $|c(E)|$ used 
on $G = (V, E)$ satisfies $|E| + |c(E)| \geq \binom{n}{2}$, then 
$(G, c)$ admits a rainbow $C_3$.   
These extensions relate to a classical anti-Ramsey\footnote{For a 
comprehensive survey of anti-Ramsey theory, see~\cite{FMO}.}   
result of Erd\H{o}s et al.~\cite{ESS} 
that any edge-coloring 
of $G = K_n$ with $n$ colors
admits a rainbow $C_3$.  
More generally, the following holds.  

\begin{theorem}[Montellano-Ballesteros and Neumann-Lara~\cite{BN}]  
\label{thm:MBNL}  
For every integer $\ell \geq 3$, 
every 
edge-colored complete graph $(K_n, c)$
satisfying 
$$
|c(E(K_n))| \geq \left(\frac{\ell-2}{2} + \frac{1}{\ell-1}\right)n + O(1)
$$
admits a rainbow $C_{\ell}$.  
\end{theorem}  

Theorem~\ref{thm:MBNL}  
confirmed a conjecture in~\cite{ESS} whose sharpness was already noted  
there:  
let $n = q(\ell - 1) + r$ for $q, r \in \mathbb{Z}$ satisfying $0 \leq r < \ell-1$;   
let $V(K_n) = V_1 \dot\cup \dots \dot\cup V_{q+1}$ be a partition 
satisfying 
$|V_1| = \dots = |V_q| = \ell-1$ and $|V_{q+1}| = r$; 
let         
all pairs of $\bigcup_{1 \leq i \leq q+1} \binom{V_i}{2}$ be given distinct colors;  
let 
all pairs crossing $V_i$ and $V_{i+1} \dot\cup \dots 
\dot\cup V_{q+1}$ be given a new color $\xi_i$, where $1 \leq i \leq q$.  
This coloring is locally imbalanced, so one may seek bounds on 
$\delta^c(K_n)$ ensuring a rainbow $C_{\ell}$ in $(K_n, c)$.  
For fixed $\ell \geq 3$, 
Axenovich et al.~\cite{AJT} 
proved that 
$\delta^c(K_n) \geq (1 + o(1))n/2$ ensures a rainbow 
$C_{\ell}$, but that for $\ell = 3$ 
the bound 
$\delta^c(K_n) \geq (1 + o(1)) \log_2 n$ 
already suffices (where $\log_2 n$
is necessary).  
Thus, 
replacing $K_n$ with an $n$-vertex host $G = (V, E)$ 
(see Theorem~1.1) 
significantly
changes the nature of the problem.  

\section{Proof of Theorem~\ref{thm:Li}}  
\label{sec:Li}  
We recall Li's proof~\cite{li2013rainbow} of Theorem~\ref{thm:Li}.
Let $(G, c)$ be an $n$-vertex 
edge-colored and edge-minimal   
graph with no rainbow
triangle $C_3$.  
We show that $\delta^c(G) \leq n/2$.  To that end,   
for a color $\alpha \in c(E)$ and a vertex $v \in V$, we define the
{\it $\alpha$-neighborhood} 
$$
N_{\alpha}(v) = \left\{ u \in N(v): \, c(\{u, v\}) = \alpha\right\}, 
$$
where $N(v) = N_G(v) = \{u \in V: \, \{u, v\} \in E\}$ 
is the usual neighborhood of $v$ in $G$, 
and $N[v] = \{v\} \cup N(v)$ is the {\it closed neighborhood} of $v$ 
in $G$.  
We define 
$$
N_! (v) = \bigcup_{\alpha \in c(E)} \left\{N_{\alpha}(v): \, 
|N_{\alpha}(v)| = 1 \right\}
$$
for the set of neighbors $u \in N(v)$ for which $c(\{u, v\})$ appears
uniquely among $\{v, w\} \in E$.  
We 
define the {\it replication number} $R = R(G, c)$ of $(G, c)$ by 
\begin{equation}
\label{eqn:replication}  
R = R(G, c) = 
\max_{v \in V} \max_{\alpha \in c(E)} 
|N_{\alpha}(v)|.      
\end{equation}  
For $v \in V$ and $U \subseteq V$, we denote
by $\deg_G^c(v, U)$ the number of colors $c(\{u, v\})$ 
among $u \in N(v) \cap U$.

Fix $(z, \zeta) \in V \times c(E)$ for which 
$|N_{\zeta}(z)| = R$ (cf.~\eqref{eqn:replication}).  
If $N_!(z) = \emptyset$, then 
each color incident to $z$ appears at least twice, so 
$\delta^c(G) \leq \deg_G^c(z) \leq (n-1)/2$ 
follows.  
Henceforth, we assume $N_!(z) \neq \emptyset$, and we define the directed graph
$D = (V_D, \vec{E}_D)$  
on vertex set $V_D = N(z)$ by putting, 
for each edge $\{x, y\} \in E(G[N(z)])$, the arc 
$(x, y) \in \vec{E}_D$ if, and only if, 
$y \in N_!(z)$ and 
$c(\{x, z\}) = c(\{x, y\})$.  
Then 
\begin{equation}
\label{eqn:7.3.2019.7:08p}  
|V_D| = 
\deg_G(z) \geq \deg_G^c(z) + R - 1, \qquad \text{and} \qquad   
|\vec{E}_D| = 
\sum_{y \in N_!(z)} \deg_D^-(y) = 
\sum_{x \in N(z)} \deg_D^+(x),     
\end{equation}  
where 
$\deg_D^+(x)$ denotes the {\sl out-degree}   
of a vertex $x \in V_D$ in $D$,  
and $\deg_D^-(x)$ denotes the corresponding {\sl in-degree}.  
We make three observations on $D$:
\smallskip 
\begin{enumerate}
\item[$(i)$]  
Each $(x, y) \in \vec{E}_D$ places $x \in N_!(z)$, lest  
some $x \neq x' \in N(z) \setminus N_!(z)$ gives 
$c(\{x', z\}) = c(\{x, z\}) = c(\{x, y\})$    
(cf.~Remark~\ref{rmk:edgeminimal});    
\smallskip 
\item[$(ii)$]  
Each $x \in N_!(z)$ 
with $\alpha = c(\{x, z\})$ 
satisfies $\deg_D^+(x) = |N_{\alpha}(x) \cap N_!(z)| \leq R - 1$
(cf.~(\ref{eqn:replication}));   
\item[$(iii)$]  
Each $y \in N_!(z)$ has 
$\deg_G^c(y, N[z]) \leq 1 + \deg_D^-(y)$, as 
$c(\{x, y\}) = \beta \neq \alpha = c(\{y, z\})$
for $x \in N(z)$ puts 
$(x, y) \in \vec{E}_D$,   
since 
$c(\{x, z\}) \neq \alpha$ by $y \in N_!(z)$ 
and $c(\{x, z\}) = \beta$
lest $\{x, y, z\}$ is rainbow.  
\end{enumerate}  
\smallskip 
\noindent Thus, 
$$
\sum_{y \in N_!(z)} (\deg_G^c(y, N[z]) - 1 ) 
\stackrel{(iii)}{\leq}  
\sum_{y \in N_!(z)} \deg_D^-(y)  
\stackrel{\eqref{eqn:7.3.2019.7:08p}}{=}    
\sum_{x \in N(z)} \deg_D^+(x) 
\stackrel{(i)}{=}  
\sum_{x \in N_!(z)} \deg_D^+(x) 
\stackrel{(ii)}{\leq}  
|N_!(z)| (R - 1)\,.    
$$
Averaging over $N_!(z)$ guarantees a vertex $y_0 \in N_!(z)$
for which 
\begin{equation}
\label{eqn:7.5.2019.2:36p}  
\deg_G^c(y_0, N[z]) \leq R.  
\end{equation}  
Since $\deg_G^c(y_0, V \setminus N[z]) \geq 
\deg_G^c(y_0) - \deg_G^c(y_0, N[z])$, we conclude    
\begin{multline*}  
n - 1 - \deg_G(z) = 
n - |N[z]|  \geq 
\deg_G^c(y_0, V \setminus N[z]) \geq 
\deg_G^c(y_0) - \deg_G^c(y_0, N[z])   
\stackrel{\eqref{eqn:7.5.2019.2:36p}}{\geq} \deg_G^c(y_0) - R  \\
\implies 
n - 1 \geq \deg_G(z) + \deg_G^c(y_0) - R 
\stackrel{\eqref{eqn:7.3.2019.7:08p}}{\geq}    
\deg_G^c(z) + R - 1 + \deg_G^c(y_0) - R 
\end{multline*}  
from which $2 \delta^c(G) \leq \deg_G^c(z) + \deg_G^c(y_0) \leq n$ 
and $\delta^c(G) \leq n/2$ follow.  \hfill $\Box$

\section{Tools for proving Theorem~\ref{thm:odd}}  
\label{sec:tools}  
All tools of this section depend on 
the following concepts of {\it separation} and {\it restriction}.   

\begin{definition}[separates/restricts]  
\label{def:restriction}  
\rm
Let $(G, c)$ be an edge-colored graph, and fix $v \in V = V(G)$ 
and $X \subseteq N(v)$.  
We say a color $\alpha \in c(E)$ will {\it $X$-separate} a vertex $y \in V$ 
from $v$ 
when 
some $x \in N(y) \cap X$ satisfies 
$\alpha = c(\{x, y\}) \neq c(\{v, x\})$. 
If, additionally, 
$\alpha \neq c(\{w, y\})$ for all $w \in N(y) \setminus X$, then we say 
that $(v, X)$ {\it restricts} the color $\alpha$ for $y$.    
We denote by $\sigma_{v,X}(y)$ the number of colors $\alpha \in c(E)$
which $X$-separate $y$ from $v$, 
and we denote by $\rho_{v,X}(y)$ 
the number of colors $\alpha \in c(E)$ 
restricted for $y$ by $(v, X)$.
\end{definition}

Every color $\alpha \in c(E)$ restricted for $y$ by $(v, X)$ also 
$X$-separates $y$ from $v$, so   
$\sigma_{v,X}(y) \geq \rho_{v,X}(y)$ holds.  
The next result 
formally extends Theorem~\ref{thm:Li}  
(see Remark~\ref{rem:Liext})    
by averaging these numbers.

\begin{prop}
\label{prop:1}  
Let $(G, c)$ be an $n$-vertex edge-colored and edge-minimal 
(cf.~Remark~\ref{rmk:edgeminimal}) 
graph
with 
$R = R(G, c)$ from~\eqref{eqn:replication}, and fix     
$v \in V$, $X \subseteq N(v)$, and 
$\emptyset \neq Y \subseteq V \setminus \{v\}$.  Then  
$$
\frac{1}{|Y|} \sum_{y \in Y} \sigma_{v,X}(y) 
\geq 
\frac{1}{|Y|} \sum_{y \in Y} \rho_{v,X}(y) 
\geq \delta^c(G) + |X| - n - (R-1)\frac{|X \cap N_!(v)|}{|Y|}.  
$$
\end{prop}

\begin{proof}[Proof of Proposition~\ref{prop:1}]
Let $(G, c)$, $R$, $v$, $X$ and $Y$ be given as above, where 
it suffices to prove the rightmost 
inequality for $X \neq \emptyset$.   
Define the directed graph $D = (V_D, \vec{E}_D)$ on vertex set 
$V_D = X \cup Y$ by putting, 
for each edge $\{x, y\} \in E$
with $x \in X$ and $y \in Y$, 
the arc 
$(x, y) \in \vec{E}_D$ 
if, and only if, 
$c(\{x, y\}) = c(\{v, x\})$.  
Similarly to~$(i)$ and~$(ii)$ of Section~\ref{sec:Li}, 
each 
$(x, y) \in \vec{E}_D$ gives $x \in N_!(v)$ 
and $\deg_D^+(x) \leq R - 1$, so 
\begin{equation}
\label{eqn:7.6.2019.1:40p}  
\sum_{y \in Y} \deg_D^-(y) = 
|\vec{E}_D| = 
\sum_{x \in X} \deg_D^+(x) = 
\sum_{x \in X \cap N_!(v)} \deg_D^+(x) 
\leq 
(R - 1) |X \cap N_!(v)|.     
\end{equation}  
Similarly to~$(iii)$ of Section~\ref{sec:Li}, 
each $y \in Y$   
admits at most 
$\deg_D^-(y) + \rho_{v,X}(y)$  
many 
colors $\alpha \in c(E)$:   
\begin{enumerate}
\item[$(a)$]  
$\alpha = c(\{x, y\})$ for some $x \in N(y) \cap X$; 
\item[$(b)$]  
$\alpha \neq c(\{w, y\})$ for all $w \in N(y) \setminus X$.  
\end{enumerate}  
Indeed, 
let 
$\alpha = c(\{x, y\})$ be such a color.  
If 
$\alpha = c(\{x, y\}) = c(\{v, x\})$, then 
$(x, y) \in \vec{E}_D$, and 
otherwise 
$(v, X)$ restricts 
$\alpha = c(\{x, y\}) \neq c(\{v, x\})$ 
for $y$  
(cf.~Definition~\ref{def:restriction}).
Consequently, 
\begin{multline}
\label{eqn:7.6.2019.1:27p}  
n - |X|
\geq 
\deg_G^c(y, V \setminus X) \geq \deg_G^c(y) - \deg_D^-(y) - \rho_{v,X}(y) \\
\implies 
\qquad 
\deg_D^-(y) \geq \deg_G^c(y) - \rho_{v, X}(y) + |X| - n
\geq \delta^c(G) - \rho_{v,X}(y)  + |X| - n.
\end{multline}  
Applying~\eqref{eqn:7.6.2019.1:27p}  
to~\eqref{eqn:7.6.2019.1:40p} renders the desired result.    
\end{proof}

\begin{remark}
\label{rem:Liext}  
\rm 
Proposition~\ref{prop:1} implies Theorem~\ref{thm:Li}:  
Let $(G, c)$ be edge-minimal with no rainbow $C_3$, 
and fix $(z, \zeta)$ with $|N_{\zeta}(z)| = R$,    
$x \in N(z) = X$, and (if possible) $y \in N_!(z) = Y$.  
Then $\rho_{z, X}(y) = 0$ 
as 
$c(\{x, y\}) \neq c(\{x, z\})$ 
gives $c(\{y, z\}) = c(\{x, y\})$ with 
$z \not\in X$,  
since $\{x, y, z\}$ is not rainbow and 
$c(\{y, z\}) = c(\{x, z\})$ violates $y \in N_!(z)$.  
Now, 
$\delta^c(G) \leq n - |X| + R - 1 \leq n - \delta^c(G)$
so $\delta^c(G) \leq n/2$.  \hfill $\Box$  
\end{remark}

For $(v, X)$ fixed, 
Proposition~\ref{prop:1} shows that some vertices $y \in V$ may admit many
colors which $X$-separate $y$ from $v$.  
For relevant $(G, c)$, 
Proposition~\ref{prop:2} finds vertices 
$y \in V$ with few such colors.  

\begin{prop}
\label{prop:2}
Fix an integer $\ell \geq 3$, 
and 
let $(G, c)$ be an edge-colored 
and edge-minimal\footnote{Recall again  
Remark~\ref{rmk:edgeminimal}.}  
graph with no rainbow $\ell$-cycle
$C_{\ell}$.  
Fix $v \in V = V(G)$ and $X \subseteq N(v)$, and let 
$\cC_{\rep} = \cC_{\rep}(v, X)$  
be the colors 
that 
repeat 
on the edges between $v$ and $X$.
Let $Y = Y(v, \cC_{\rep})$ be the vertices $y \in V$ 
to which there is an $(\ell-1)$-vertex rainbow path $P_{vy}$ 
from $v$, none of whose edges has a color in $\cC_{\rep}$.   
Then every $y\in Y$ 
satisfies $\sigma_{v,X}(y) \leq 3\ell$.  
\end{prop}

\begin{proof}[Proof of Proposition~\ref{prop:2}]  
Let $(G, c)$, $v$, $X$, $\cC_{\rep}$, 
$y \in Y$, 
and 
$P_{vy}$ be given as above.   
For a vertex $x \in N(y) \cap X$,  
the subgraph $P_{vy} + \{x, y\} + \{v, x\}$ is a rainbow $\ell$-cycle in $(G, c)$ 
unless:  
$$
\text{$(A)$ 
$x \in V(P_{vy})$;} \enskip
\text{$(B)$ $c(\{x, y\}) \in c(E(P_{vy}))$;} \enskip  
\text{$(C)$ $c(\{v, x\}) \in c(E(P_{vy}))$;} \enskip \text{or} \enskip
\text{$(D)$ $c(\{x, y\}) = c(\{v, x\})$.}    
$$
At most $|N(y) \cap X \cap V(P_{vy})| \leq \ell - 3$ 
colors $c(\{x, y\})$ 
are given by a vertex $x \in N(y) \cap X \cap V(P_{vy})$ 
satisfying~$(A)$ and at most $|E(P_{vy})| \leq \ell - 2$ colors 
$c(\{x, y\})$ satisfy~$(B)$.    
At most $|E(P_{vy})| \leq \ell - 2$ colors $c(\{x, y\})$ satisfy~$(C)$
because  $c(\{v, x\}) \not\in \cC_{\rep}$.  
All remaining $c(\{x, y\})$ over $x \in N(y) \cap X$ satisfy~$(D)$, 
lest $(G, c)$ admits a rainbow $\ell$-cycle $C_{\ell}$.  
\end{proof}

\subsection{Some corollaries}  
\label{sec:corollaries}  
We now consider several useful corollaries of Propositions~\ref{prop:1}
and~\ref{prop:2}.  

\begin{cor}
\label{cor:3}  
Fix an integer $\ell \geq 3$,  
and 
let $(G, c)$ be an $n$-vertex edge-colored and edge-minimal graph 
with 
no 
rainbow
$\ell$-cycle $C_{\ell}$.  
Let $(z, \zeta) \in V \times c(E)$ satisfy 
$|N_{\zeta}(z)| = R$ (cf.~\eqref{eqn:replication}), and     
let $Y = Y(z,\zeta) \subseteq V$ be the vertices
$y \in V$ 
 to which there is an $(\ell-1)$-vertex rainbow path $P_{zy}$ 
from $z$, none of whose edges is colored $\zeta$.  
If $Y \neq \emptyset$, 
$$
\delta^c(G) \leq 
\frac{n}{2} + 
\max 
\left\{0, \, 
3\ell + 
(R - 1)\left(\frac{n+1}{2|Y|} - 1\right)    
\right\}.  
$$
\end{cor}  

\begin{proof}[Proof of Corollary~\ref{cor:3}]   
Let $(G, c)$, $z$, $\zeta$, $R$ and $Y = Y(z,\zeta) \neq \emptyset$ 
be given as above, where for the sake of an argument we assume $\delta^c(G) > n/2$.  
Let $X \subseteq N(z)$ satisfy that $|X| = \lceil n/2 \rceil$, 
that 
$\zeta = c(\{x_0, z\})$ for some $x_0 \in X$, and that 
all $\{x, z\}$ with $x \in X$ are colored distinctly.   
Set $X^+ = X \cup N_{\zeta}(z)$,  
and set $\cC_{\rep} = \cC_{\rep}(z, X^+)$ to be the colors $\alpha = c(\{x, z\})$
repeating among $x \in X^+$.  
Then $\cC_{\rep} \subseteq \{\zeta\}$, which  
by hypothesis 
is forbidden
on the path   
$P_{zy}$ ending in $y \in Y = Y(z,\zeta)$.  
Proposition~\ref{prop:2} guarantees that $\sigma_{z, X^+}(y) \leq 3\ell$  
holds 
for each $y \in Y$, and   
Proposition~\ref{prop:1} then renders
$$
3\ell
\geq 
\frac{1}{|Y|}  
\sum_{y \in Y} \sigma_{z,X^+}(y) 
\geq \delta^c(G) + |X^+| - n - \frac{(R-1) |X^+ \cap N_!(z)|}{|Y|},       
$$
and 
using $|X^+| = |X| + R - 1$ and $\lceil n/2 \rceil = |X| \geq |N_!(z) \cap X^+|$ 
completes the proof.  
\end{proof}

In practice, the set $Y = Y(z, \zeta)$ 
in Corollary~\ref{cor:3}  
will be large, 
and will guarantee
the following result.

\begin{cor}
\label{cor:4}  
Fix an integer $\ell \geq 3$, 
and 
let $(G, c)$ be an $n$-vertex edge-colored graph 
with no rainbow $\ell$-cycle $C_{\ell}$.  
Then 
$\delta^c(G) \leq (n/2) + 3\ell$.  
\end{cor}  

\begin{proof}[Proof of Corollary~\ref{cor:4}]
Let $(G, c)$ be given as above.
For the sake of an argument, we assume $\delta^c(G) \geq (n/2) + 2\ell - 5$,    
and w.l.o.g.~we assume $(G, c)$ is edge-minimal.    
Let $(z, \zeta) \in V \times c(E)$ satisfy that $|N_{\zeta}(z)| = R$ (cf.~\eqref{eqn:replication}).    
For $1 \leq i \leq \ell - 1$, 
let $Y_i = Y_i(z,\zeta)$ be the set of vertices $y_i \in V$ 
to 
which 
there is an $i$-vertex rainbow path $P_{zy_i}$ from $z$, none of whose edges
is colored $\zeta$.  
Inductively, these sets are non-empty as 
$Y_1 = \{z\}$, 
and for some 
$1 \leq j \leq \ell - 2$, a fixed   
$y_j \in Y_j$
and corresponding path $P_{z y_j}$   
provide 
\begin{equation}
\label{eqn:7.7.2019.5:19p}  
|Y_{j+1}| \geq \deg_G^c(y_j) - 1 - |E(P_{zy_j})| - (|V(P_{zy_j})| - 2)
\geq \delta^c(G) - 2j + 2  \geq \delta^c(G) - 2\ell + 6 \geq \tfrac{n+1}{2}
\end{equation}  
with $\delta^c(G) \geq (n/2) + 2\ell - 5$.  
Corollary~\ref{cor:3} now guarantees 
$$
\delta^c(G) \leq \frac{n}{2} + 3\ell + (R-1)
\left(\frac{n+1}{2|Y_{\ell-1}|} - 1 \right)
\stackrel{\eqref{eqn:7.7.2019.5:19p}}{\leq}    
\frac{n}{2} + 3\ell, 
$$
as desired.  
\end{proof}

The following corollary describes sets similar to $Y(z, \zeta)$
which are also large.

\begin{cor}
\label{cor:5}
Fix an integer $\ell \geq 3$, 
and 
let $(G, c)$ be an edge-colored and edge-minimal graph
with no rainbow $\ell$-cycle $C_{\ell}$.  
Let $T$ be a triangle in $G$, let $v \in V(T)$, and let 
${\cal C}_T \subseteq c(E)$ satisfy
$\cC_T \cap c(E(T)) = \emptyset$.    
Let $Y = Y(v, \cC_T) = Y_{\ell-1}(v, \cC_T)$ be the vertices $y \in Y$ 
to which there is an $(\ell-1)$-vertex 
rainbow path $P_{v y}$ from $v$, none of whose edges has a color
in $\cC_T$.  
Then $|Y| \geq (3/2)(\delta^c(G) - |\cC_T| - 4\ell)$.  
\end{cor}

\begin{proof}[Proof of Corollary~\ref{cor:5}]  
Let $(G, c)$, $T$, $v$ and ${\cal C}_T$ be given 
as above, 
where for sake of argument we assume $\delta^c(G) \geq |\cC_T| + 4\ell + 1$, 
and      
where we set $\hat{\cC}_T = \cC_T \cup c(E(T))$.  
Since $T$ is not monochromatic 
(cf.~Remark~\ref{rmk:edgeminimal}), 
we label $V(T) = \{v, x_1, x_2\}$
with $c(\{v, x_2\}) \neq c(\{x_1, x_2\})$.  
For $1 \leq i \leq \ell - 1$, 
let 
$W_i = W_i(x_1, \hat{\cC}_T)$ 
be the set of vertices $w_i \in V$ 
to which there is an $i$-vertex rainbow path $P_{x_1 w_i}$ 
from $x_1$, none of whose edges has a color in $\hat{\cC}_T$, and whose
vertices meet $V(T)$ only in $x_1$.    
Inductively, these sets are non-empty as 
$W_1 = \{x_1\}$, 
and for some $1 \leq j \leq \ell-2$,  
a fixed $w_j \in W_j$
and corresponding path $P_{x_1 w_j}$  
provide 
that 
\begin{equation}
\label{eqn:7.9.2019.1:24p}  
|W_{j+1}| \geq \deg_G^c(w_j) - |\hat{\cC}_T|
 - |E(P_{x_1 w_j})| 
- (|V(P_{x_1 w_j})| - 2) - |\{v, x_2\}| 
\geq \deg_G^c(w_j) - 2(j+1) - |\cC_T|   
\end{equation}  
is positive from $\delta^c(G) \geq |\cC_T| + 4\ell + 1$.  
It is easy to see that 
\begin{equation}
\label{eqn:7.9.2019.2:30p}  
W_{\ell-3}(x_1, \hat{\cC}_T) 
\cup 
W_{\ell-2}(x_1, \hat{\cC}_T) 
= 
W_{\ell - 3} \cup 
W_{\ell-2} 
\subseteq Y = Y_{\ell-1}(v, \cC_T).   
\end{equation}  
Indeed, if $w_{\ell-3} \in W_{\ell - 3}$ is given by 
$P_{x_1 w_{\ell-3}}$, 
then 
the path 
$P_{x_1 w_{\ell-3}} + \{v, x_2\} + \{x_1, x_2\}$ places 
$w_{\ell-3} \in Y$, and if $w_{\ell-2} \in W_{\ell-2}$ is given by 
$P_{x_1 w_{\ell-2}}$, then the path   
$P_{x_1 w_{\ell-2}} + \{v, x_1\}$ places 
$w_{\ell-2} \in Y$.    
We bound~\eqref{eqn:7.9.2019.2:30p}  
as follows.   
Let $\Gamma = G[W_{\ell-3}]$ be the edge-colored subgraph of $G$
induced on $W_{\ell-3}$.  
Then $\Gamma$ admits no rainbow $\ell$-cycles $C_{\ell}$, whence 
Corollary~\ref{cor:4}
guarantees a vertex
$w_{\ell-3} \in W_{\ell-3}$ for which $\deg_{\Gamma}^c(w_{\ell-3}) \leq (1/2)
|W_{\ell-3}| + 3\ell$.  
As such 
(see the last inequality of~\eqref{eqn:7.9.2019.1:24p} to bound $|W_{\ell-2}|$), 
\begin{equation}  
\label{eqn:7.9.2:33p}  
|W_{\ell-2} \setminus W_{\ell-3}| 
\geq 
\deg_G^c(w_{\ell-3}) 
- 2(\ell-2) 
- |\cC_T|  
- \deg_{\Gamma}^c(w_{\ell-3})
\geq 
\delta^c(G) - \tfrac{1}{2} |W_{\ell-3}| - 5\ell - |\cC_T|, 
\end{equation}  
and so 
$$
|Y| 
\stackrel{\eqref{eqn:7.9.2019.2:30p}}{\geq}  
|W_{\ell-2} \cup W_{\ell-3}|  = |W_{\ell-2} \setminus W_{\ell-3}| + |W_{\ell-3}|
\stackrel{\eqref{eqn:7.9.2:33p}}{\geq}    
\delta^c(G) + \tfrac{1}{2} 
|W_{\ell-3}|  
- 5\ell - |\cC_T|  
\stackrel{\eqref{eqn:7.9.2019.1:24p}}{\geq}    
\tfrac{3}{2} \delta^c(G) - 6 \ell - \tfrac{3}{2} |\cC_T|,   
$$
as promised.  
\end{proof}

\section{Proof of Theorem~\ref{thm:odd}}  
\label{sec:odd}  
Fix an integer\footnote{By Theorem~\ref{thm:even}, it suffices to prove
Theorem~\ref{thm:odd} for odd integers $\ell$.  However, most of the current 
argument
is independent of parity considerations, so we make
no distinction now.}  
$\ell \geq 3$.  
Let $(G, c)$ be an $n$-vertex  
edge-colored and edge-minimal graph 
(cf.~Remark~\ref{rmk:edgeminimal})   
satisfying $\delta^c(G) \geq (n+1)/2$.    
We assume that $(G, c)$ admits no rainbow $\ell$-cycle
$C_{\ell}$, and we bound $n \leq n_0(\ell)$ from above in the course of this proof.
Fix $(z, \zeta) \in V \times c(E)$ with $|N_{\zeta}(z)| = R$ (cf.~\eqref{eqn:replication}).    
Let $X \subset N(z)$ satisfy that $|X| = \delta^c(G) - 1$ and that 
$c(\{x, z\}) \neq \zeta$ are distinct among $x \in X$.   
We distinguish two cases.

\medskip

\noindent {\bf Case 1 ($\exists\,  e_0 \in E(G[X]): \, c(e_0) \neq 
\zeta$).} 
By our choice of $X$, 
the following hold:  
\begin{enumerate}
\item[(I)]  
$\zeta$ does not appear on 
the triangle $T = \{z\} \cup e_0$; 
\item[(II)]  
$\zeta$ 
is 
the only color 
possibly 
repeating among $c(\{x,z\})$ for $x \in X^+ = X \cup N_{\zeta}(z)$.
\end{enumerate}  
As such, we set $\cC_{\rep} = \cC_T \subseteq 
\{\zeta\}$ so that 
the set 
$Y = 
Y(z, \cC_{\rep}) = 
Y(z, \zeta) = Y(z, \cC_T)$
commonly featured in each of Proposition~\ref{prop:2} and Corollaries~\ref{cor:3}
and~\ref{cor:5} has size, by the last of these,  
\begin{equation}
\label{eqn:7.15.2019.5:13p}  
|Y| \geq \tfrac{3}{2} (\delta^c(G) - 1 - 4\ell)
\geq \tfrac{3}{2} \left(\tfrac{n+1}{2} - 1 - 4\ell\right) 
\stackrel{*}{\geq} \tfrac{2}{3} 
(n+1), 
\end{equation}    
where $*$ holds 
when $n \geq 78\ell$, which we assume for the sake of an argument.  
Corollary~\ref{cor:3} then yields 
\begin{equation}
\label{eqn:7.16.2019.12:18p}  
\frac{n+1}{2} \leq \delta^c(G) \leq \frac{n}{2} + 3\ell + 
(R - 1)\left(\frac{n+1}{2|Y|} - 1\right)  
\stackrel{\eqref{eqn:7.15.2019.5:13p}}{\leq}    
\frac{n}{2} + 3\ell - \frac{1}{4}(R-1)
\quad \implies \quad R \leq 12\ell.  
\end{equation}

Now, define the directed graph $F = (V, \vec{E}_F)$ 
on vertex set $V = V(G)$, 
where 
$$
\vec{E}_F = \left\{(x, y) \in X^+ \times V: \, 
\{x, y\} \in E = E(G) \text{ and } 
c(\{x, y\}) \neq c(\{x, z\})\right\}.  
$$
Note that for every $(x, y) \in \vec{E}_F$, the color $c(\{x, y\})$ does $X^+$-separate $y$ from $z$  
(cf.~Definition~\ref{def:restriction}).    
On the one hand, 
every $x \in X^+$ clearly satisfies $\deg_F^+(x) \geq \deg_G^c(x) - 1$, and so 
\begin{equation}  
\label{eqn:7.15.2019.3:59p}  
|\vec{E}_F| = 
\sum_{x \in X^+} \deg_F^+(x) 
\geq \sum_{x \in X^+} (\deg_G^c(x) - 1)
\geq |X^+| (\delta^c(G) - 1).  
\end{equation}  
On the other hand, 
with 
$Y$ 
defined above, 
\begin{equation}
\label{eqn:7.19.2019.2:38p}  
|\vec{E}_F| = \sum_{y \in V} \deg_F^-(y) 
= 
\sum_{y \in V \setminus Y} \deg_F^-(y)  
+ 
\sum_{y \in Y} \deg_F^-(y)  
\leq 
(n - |Y|) |X^+| + \sum_{y \in Y} 
\deg_F^-(y).  
\end{equation}  
For a fixed $y \in Y$, 
we bound 
$$
\deg_F^-(y)  =  \left|\left\{x \in N(y) \cap X^+:  \, 
c(\{x, y\}) \neq c(\{x, z\}) \right\}\right|  
 =  
\sum_{\alpha \in c(E)} 
\big|  
\left\{x \in N_{\alpha}(y) \cap X^+: \,  c(\{x, z\}) \neq \alpha  
\right\}
\big|.  
$$
Let $\cA = \cA_y$ be these colors $\alpha \in c(E)$ admitting some 
$x \in N_{\alpha}(y) \cap X^+$ with 
$c(\{x, z\}) \neq \alpha$ (where $\alpha = c(\{x, y\})$ from $x \in N_{\alpha}(y)$).  Then $\cA$ is precisely the set of colors which $X^+$-separate $y$ from $z$, so 
$|\cA| = \sigma_{z, X^+}(y)$ 
holds 
by Definition~\ref{def:restriction}.    
Then 
\begin{equation}
\label{eqn:7.19.2019.1:37p}  
\deg_F^-(y)  
 =  
\sum_{\alpha \in \cA} 
\big|  
\left\{x \in N_{\alpha}(y) \cap X^+: \, \alpha \neq c(\{x, z\})
\right\}
\big|  
\leq 
\sum_{\alpha \in \cA} 
\big|  
N_{\alpha}(y) \cap X^+\big|
\leq 
\sum_{\alpha \in \cA} 
|N_{\alpha}(y)\big|  
  \stackrel{\eqref{eqn:replication}}{\leq}
|\cA|R.  
\end{equation}  
Since $y \in Y = Y(z, \cC_{\rep})$, Proposition~\ref{prop:2}
guarantees that $|\cA| = \sigma_{z, X^+}(y) \leq 
3\ell$, and so 
\begin{equation}
\label{eqn:7.21.2019.1:36p}  
\deg_F^-(y)   
\stackrel{\eqref{eqn:7.19.2019.1:37p}}{\leq}    
|\cA|R
= \sigma_{z, X^+}(y) \cdot R
\stackrel{\text{Prop.\ref{prop:2}}}{\leq}  
3\ell R
\stackrel{\eqref{eqn:7.16.2019.12:18p}}{\leq}  
36 \ell^2.
\end{equation}  
Applying~\eqref{eqn:7.21.2019.1:36p}  
to~\eqref{eqn:7.19.2019.2:38p}  
yields 
\begin{equation}
\label{eqn:7.16.2019.12:42p}  
|\vec{E}_F| 
\leq 
(n - |Y|) |X^+| + \sum_{y \in Y} 
\deg_F^-(y)  
\leq 
(n - |Y|) |X^+| + 36\ell^2 |Y|.   
\end{equation}  
Comparing~\eqref{eqn:7.15.2019.3:59p}  
and~\eqref{eqn:7.16.2019.12:42p}  
yields 
$|X^+| (\delta^c(G) - 1) \leq (n - |Y|) |X^+| + 36\ell^2 |Y|$, 
or equivalently, 
$$
n \geq 
\delta^c(G) - 1 
+ 
\left(1 - \frac{36\ell^2}{|X^+|}\right) |Y|.       
$$
Using $|X^+| = \delta^c(G) - 1 + R \geq \delta^c(G) \geq (n+1)/2$, we infer 
\begin{equation}
\label{eqn:7.19.2019.2:58p}    
\tfrac{1}{2} (n + 1)  
\geq n - \delta^c(G) + 1   
\geq \left(1 - \tfrac{72\ell^2}{n+1}\right) |Y| 
\stackrel{\eqref{eqn:7.15.2019.5:13p}}{\geq}    
\left(1 - \tfrac{72\ell^2}{n+1}\right)
\times 
\tfrac{2}{3}(n+1),   
\end{equation}  
which implies $n \leq 288 \ell^2 - 1$.  \hfill $\Box$   

\medskip

\noindent {\bf Case 2 ($\forall e \in E(G[X]), \, c(e) = \zeta$).}  
Set 
$Y = V \setminus (\{z\} \cup X)$.  We first observe 
\begin{equation}
\label{eqn:7.10.2019.4:26}
\delta^c(G) - 2 \leq 
|Y|
\leq \delta^c(G) - 1
\qquad \text{and} \qquad 
\tfrac{1}{2} (n-1) \leq \delta^c(G) - 1 = |X| \leq \tfrac{1}{2}(n+1).  
\end{equation}  
Indeed, 
$|Y| = n - 1 - |X| = n - \delta^c(G) \leq \delta^c(G) - 1$ holds from 
$|X| = \delta^c(G) - 1 \geq (1/2)(n-1)$.  
Now, fix $x \in X$ and $\{x, y\} \in E$ where 
$c(\{x, y\}) \neq \zeta$ and $c(\{x, y\}) \neq c(\{x, z\})$.
Then $y \in Y$ and there are 
at least $\deg_G^c(x) - 2$ many such edges.    
Thus,   
$|Y| \geq \deg_G^c(x) - 2 \geq \delta^c(G) - 2$ holds and      
$|X| \leq (1/2)(n+1)$  follows.

We now define two subsets of $Y$ that we wish to later avoid.  
For that, let $H = G[X, Y]$ be the bipartite subgraph of $G$ 
induced by the bipartition 
$X \cup Y$, and let    
$D$ be the subgraph of $H$ consisting of edges $\{x, y\} \in E(H)$ 
with 
$x \in X$, $y \in Y$, and 
$c(\{x, y\}) = c(\{x, z\})$.  
Let $Y_{H}$ be the vertices $y \in Y$
sending 
$\deg_H^c(y) \leq (5/2) \ell$  
many 
distinct colors to $X$, and    
let $Y_{D}$ be the vertices $y \in Y$ sending 
$\deg_D(y) \geq 2$ many $D$-edges to $X$.   

\begin{claim}
\label{clm:7.10.2019.6:36p}  
$|Y_{H}| \leq 11 \ell$ 
and 
$|Y_{D}| \leq |X|/2$.    
\end{claim}

\begin{proof}[Proof of Claim~\ref{clm:7.10.2019.6:36p}]  
Let $\Gamma = G[A]$
be the edge-colored subgraph of $G$ 
induced on $A = Y_{H}$.   
Since $G[A]$ has no rainbow $\ell$-cycles $C_{\ell}$, 
Corollary~\ref{cor:4} guarantees $a \in A$ with 
$\deg_{\Gamma}^c(a) \leq (1/2)|A| + 3\ell$.  
Since $a$ sends at most $(5/2) \ell + 1$ distinct colors 
to $\{z\} \cup X$ and at most $|Y| - |A|$ distinct colors
to $Y \setminus A$, we see 
$$
\tfrac{1}{2} |A| + 3\ell \geq 
\deg_{\Gamma}^c(a) 
\geq 
\deg_G^c(a) - \tfrac{5}{2} \ell - 1 - |Y| + 
|A|  
\stackrel{\eqref{eqn:7.10.2019.4:26}}{\geq}
|A| 
- \tfrac{5}{2} \ell     
\qquad \implies \qquad |Y_H| = |A| \leq 11\ell.  
$$
Since 
each $x \in X$ sends to $Y$ precisely $\deg_D(x)$ many 
$D$-edges
and $\geq \deg_G^c(x) - 2$ many $\zeta$-free $H\setminus D$-edges,  
\begin{multline*}  
\deg_D(x) + \deg_G^c(x) - 2 \leq \deg_H(x) \leq |Y| 
\stackrel{\eqref{eqn:7.10.2019.4:26}}{\leq} \delta^c(G) - 1
\quad \implies \quad 
\deg_D(x) \leq 1  \\
\implies \quad 
 2|Y_D| \leq 
\sum_{y \in Y_D} \deg_D(y) \leq    
\sum_{y \in Y} \deg_D(y) = 
|E(D)| = \sum_{x \in X} \deg_D(x) \leq |X|,     
\end{multline*}  
and so $|Y_{D}| \leq |X| / 2$ follows.  
\end{proof}

Continuing with Case~2,   
set $Y_0 = Y \setminus (Y_{H} \cup Y_{D})$, 
set 
$H[X, Y_0] = G[X, Y_0]$ to be the bipartite 
subgraph of $H$ induced by the bipartition $X \cup Y_0$, 
and 
set 
$H_0 = H[X, Y_0] \setminus D$.
For each $x \in X$, we already observed (cf.~\eqref{eqn:7.10.2019.4:26}) 
that $x$ sends at least $\deg^c_G(x) - 2$
many 
non-$\zeta$, non-$c(\{x, z\})$ 
colors into $Y$, and so   
\begin{equation}
\label{eqn:7.11.2019.1:28p}  
\forall \, x \in X, \qquad 
\deg_{H_0}^c(x) \geq 
\deg_G^c(x) - 2 - |Y \setminus Y_0| 
\stackrel{\text{Clm.\ref{clm:7.10.2019.6:36p}}}{\geq}    
\delta^c(G) - 2 - 11\ell - \tfrac{1}{2}|X|  
\stackrel{\eqref{eqn:7.10.2019.4:26}}{\geq}
\tfrac{1}{4}(n+1) - 11\ell - 2.    
\end{equation}  
To $X$, 
each $y \in Y_0$ 
sends $\deg_{H}^c(y) \geq (5/2) \ell + 1$ many colors 
and 
$\deg_D(y) \leq 1$ many $D$-edges, 
and so 
\begin{equation}
\label{eqn:7.11.2019.1:29p}  
\forall \, y \in Y_0, \quad 
\deg_{H_0}^c(y) \geq \tfrac{5}{2} \ell.  
\end{equation}  

To conclude Case~2, it is convenient to now distinguish between $\ell$ (mod 2).  

\bigskip 

\noindent {\bf Case 2A ({\rm $\ell$ is odd}).}  
  With $(z, \zeta)$ fixed at the start, 
fix $y_0 \in N_{\zeta}(z)$ arbitrarily, where necessarily $y_0 \in Y$.  
The number of non-$\zeta$ colors that $y_0$ sends to $X$ is at least 
$\deg_G^c(y_0) - 1 - (|Y| - 1) \geq \delta^c(G) - |Y| 
\geq 1$ 
by~\eqref{eqn:7.10.2019.4:26}, 
so fix $x_1 \in X \cap N(y_0)$ to satisfy 
$c(\{x_1, y_0\}) \neq \zeta$.
For an even integer $k \geq 2$, let 
$Q_{k-1} = (z, y_0, x_1, y_2, \dots, x_{k-1})$
be a rainbow path, 
where $x_1, \dots, x_{k-1} \in X$
and $y_2, \dots, y_{k-2} \in Y_0$.  
Then $Q_{k-1}$ would be  extended to a rainbow path 
$Q_k = (z, y_0, \dots, x_{k-1}, y_k)$ along 
at least  
\begin{equation}
\label{eqn:7.11.2019.5:45p}  
\deg^c_{H_0}(x_{k-1}) - |E(Q_{k-1})| - |\{y_0, \dots, y_{k-2}\}| = 
\deg_{H_0}^c(x_{k-1}) - k - \tfrac{k}{2}   
\stackrel{\eqref{eqn:7.11.2019.1:28p}}{\geq}    
\tfrac{1}{4}(n+1) - 11\ell - 
2 - \tfrac{3}{2} k 
\end{equation}  
many $y_k \in Y_0 \setminus \{y_0, \dots, y_{k-2}\}$, and    
$Q_k$ would be extended to a rainbow path
$Q_{k+1} = (z, y_0, \dots, y_k, x_{k+1})$  
along at least  
$$
\deg_{H_0}^c(y_k) - |E(Q_k)| - |\{x_1, \dots, x_{k-1}\}| = 
\deg_{H_0}^c(y_k) - (k+1) -  
(k/2)
$$
many $x_{k+1} \in X \setminus \{x_1, \dots, x_{k-1}\}$.   
More strongly, $X$ was chosen with $c(\{x, z\})$ 
distinct among $x \in X$, 
so 
$Q_k$ would be extended to a rainbow path 
$Q_{k+1} = (z, y_0, \dots, y_k, x_{k+1})$ 
along at least  
\begin{equation}
\label{eqn:7.11.2019.5:46p}  
\deg_{H_0}^c(y_k) - 2(k+1) -   
(k/2) 
\stackrel{\eqref{eqn:7.11.2019.1:29p}}{\geq}    
\tfrac{5}{2} (\ell - k -  (4/5)) 
\end{equation}  
many $x_{k+1} \in X \setminus \{x_1, \dots, x_{k-1}\}$ 
where additionally  $c(\{x_{k+1}, z\}) 
 \not\in c(E(Q_k))$.  
Then $Q_{k+1}$ bears the rainbow $(k+3)$-cycle
$(z, y_0, x_1, \dots, y_k, x_{k+1}, z)$  
since $c(\{x_{k+1}, z\}) \not\in c(E(Q_k))$ holds and 
since 
$c(\{x_{k+1}, y_k\}) \neq c(\{x_{k+1}, z\})$ holds from
$\{x_{k+1}, y_k\} \not\in E(D)$.  
Since $(G, c)$ has no rainbow $\ell$-cycles $C_{\ell}$, it must be  that $k + 3 \leq \ell - 1$.  
Since~\eqref{eqn:7.11.2019.5:46p}  
is positive with $k = \ell - 4$,   
\eqref{eqn:7.11.2019.5:45p}  
must be non-positive, 
whence $n \leq 50 \ell$.  \hfill $\Box$  

\medskip

\noindent {\bf Case 2B ({\rm $\ell$ is even}).}  
The argument above slightly simplifies.
Choose $x_1 \in X$ arbitrarily.  
As before, we extend 
a rainbow path 
$\hat{Q}_{k-1} = (z, x_1, y_2, \dots, x_{k-1})$ 
with 
$x_1, \dots, x_{k-1} \in X$ and $y_2, \dots, y_{k-2} \in Y_0$  
to rainbow paths $\hat{Q}_k = (z, x_1, \dots, x_{k-1}, y_k)$
and $\hat{Q}_{k+1} = (z, x_1, \dots, y_k, x_{k+1})$ where 
$y_k \in Y_0 \setminus \{y_2, \dots, y_{k-2}\}$ and 
$x_{k+1} \in X \setminus \{x_1, \dots, x_{k-1}\}$, and where 
$c(\{x_{k+1}, z\}) \not\in c(E(\hat{Q}_k))$.  
The paths $\hat{Q}_k$ and $\hat{Q}_{k+1}$ are respectively shorter than 
$Q_k$ and $Q_{k+1}$ above, so 
inequalities analogous to those in~\eqref{eqn:7.11.2019.5:45p}  
and~\eqref{eqn:7.11.2019.5:46p}  
still hold, and with $k + 1 < \ell - 1$, 
we similarly conclude $n \leq 50 \ell$.  \hfill $\Box$

\section{Proof of Theorem~\ref{thm:ell}}  
\label{sec:ell}  
Our proof of Theorem~\ref{thm:ell} 
follows that of Theorem~\ref{thm:odd}, where we also use 
the following corollary of Propositions~\ref{prop:1}
and~\ref{prop:2} 
from Section~\ref{sec:tools}.    

\begin{cor}
\label{cor:6}
Fix an integer $\ell \geq 3$, and 
let $(G, c)$ be an $n$-vertex edge-colored and edge-minimal graph
with no rainbow $\ell$-cycle $C_{\ell}$ and 
with $\delta^c(G) \geq 5R + 27 \ell$ (cf.~\eqref{eqn:replication}).  
Then $\delta^c(G) < n/2$ or 
$\Delta(G) < \delta^c(G) + 4R + 3\ell$.  
\end{cor}  

\begin{proof}[Proof of Corollary~\ref{cor:6}]
Let $(G, c)$ be given as above.  Assume for a contradiction 
that  
$\delta^c(G) \geq n/2$ and  
that 
some 
$v \in V = V(G)$ has $\deg_G(v) \geq \delta^c(G) + 4R + 3\ell$.     
Then $\deg_G(v) \geq (n+1)/2$ whence $v$ is incident to 
some triangle $T = T_v$ in $G$, where we set\footnote{The colors $\alpha$, $\beta$, 
$\gamma$ 
are not identical by Remark~\ref{rmk:edgeminimal}, but they 
need not all be distinct.  These considerations, however, play no role in the
current context.}    
$c(E(T)) = \{\alpha, \beta, \gamma\}$.  
Since 
$|N_{\alpha}(v) \cup N_{\beta}(v) \cup N_{\gamma}(v)| \leq 3R$, 
some 
$X \subseteq N(v)$ has size $|X| = \delta^c(G) + R + 3\ell$, where
at least $\delta^c(G)$ of the colors $c(\{v, x\})$ 
are distinct among $x \in X$, 
and where 
$|N_{\alpha}(v) \cap X|, |N_{\beta}(v) \cap X|, |N_{\gamma}(v) \cap X| \leq 1$.  
Let $\cC_T = \cC_{\rep}$ be the $\leq R + 3\ell$ colors $c(\{v, x\})$
repeating among $x \in X$, where   
$\cC_T \cap c(E(T)) = \emptyset$.      
Let $Y = Y(v, \cC_{\rep}) = Y(v, \cC_T)$ be the set commonly 
featured in 
each of Proposition~\ref{prop:2} and Corollary~\ref{cor:5}.  
Corollary~\ref{cor:5} guarantees 
\begin{multline}  
\label{eqn:7.14.2019.5:14p}  
|Y| 
\geq 
\tfrac{3}{2} (\delta^c(G) - |\cC_T| - 4\ell)
\geq 
\tfrac{3}{2} (\delta^c(G) - R - 7\ell)
= 
\delta^c(G) + \tfrac{1}{2} \delta^c(G) - \tfrac{3}{2} (R + 7\ell)   \\
\geq 
\delta^c(G) + \tfrac{1}{2} (5R + 27\ell) - \tfrac{3}{2} (R + 7\ell)  
= 
\delta^c(G) + R + 3\ell = |X|.     
\end{multline}  
Proposition~\ref{prop:1} then guarantees 
\begin{multline}  
\label{eqn:7.14.2019.6:09p}  
\frac{1}{|Y|} \sum_{y \in Y} \sigma_{v,X}(y)
\geq \delta^c(G) + |X| - n - (R-1) \frac{|X \cap N_!(v)|}{|Y|}  \\
\stackrel{\eqref{eqn:7.14.2019.5:14p}}{\geq}    
\delta^c(G) + |X| - n - (R-1) 
= 2\delta^c(G) + R + 3\ell - n - (R - 1)
\geq 3\ell + 1, 
\end{multline}  
which with 
$Y = Y(v, \cC_{\rep})$ 
contradicts Proposition~\ref{prop:2}.   
\end{proof}

\subsection{Proof of Theorem~\ref{thm:ell}}  
Let $(G, c)$, $(z, \zeta)$, and $X \subset N(z)$
be given as in Section~\ref{sec:odd}.  
In Case~2, we proved that $n \leq 50 \ell$, but in Case~1 
we proved only that $n \leq 288 \ell^2$, which we now improve to $n \leq 432 \ell - 1$.
The bottleneck of Case~1 arises
in~\eqref{eqn:7.21.2019.1:36p},   
where a fixed
$y \in Y$ 
satisfies $\deg_F^-(y) \leq 3\ell R \leq 36 \ell^2$  
(cf.~\eqref{eqn:7.16.2019.12:18p}).    
We claim 
that 
\begin{equation}
\label{eqn:7.18.2019.6:28p}  
\deg_F^-(y) \leq 4R + 6\ell 
\stackrel{\eqref{eqn:7.16.2019.12:18p}}{\leq}    
54 \ell, 
\end{equation}  
which 
if true 
updates~\eqref{eqn:7.19.2019.2:58p}  
to say 
$$  
\tfrac{1}{2}(n+1) \geq \left(1 - \tfrac{108\ell}{n+1}\right)   
\times \tfrac{2}{3}(n+1),   
$$
which gives $n \leq 432 \ell - 1$.   
To see~\eqref{eqn:7.18.2019.6:28p},   
recall from~\eqref{eqn:7.19.2019.1:37p} that 
\begin{equation}
\label{eqn:7.19.2019.3:19p}  
\deg_F^-(y) \leq \sum_{\alpha \in \cA} |N_{\alpha}(y)|, 
\end{equation}  
where $\cA = \cA_y$ is the set of 
colors $\alpha \in c(E)$ where
some $x \in X^+$ satisfies $\alpha = c(\{x, y\}) \neq c(\{x, z\})$.   
In particular, $\cA$ is precisely the set of colors which $X^+$-separate $y$ from $z$, so 
$|\cA| = \sigma_{z, X^+}(y)$ 
holds 
by Definition~\ref{def:restriction}.
Moreover, 
recall (cf.~\eqref{eqn:7.21.2019.1:36p}) that   
Proposition~\ref{prop:2} guarantees $|\cA| = \sigma_{z, X^+}(y) \leq 3\ell$.       
Now, let $\cB = \cB_y$ consist of all non-$\cA$ colors 
incident to $y$, in which case 
\begin{equation}  
\label{eqn:7.18.2019.6:46p}  
\sum_{\beta \in \cB} |N_{\beta}(y)| 
\geq 
|\cB| = 
\deg_G^c(y) - |\cA|  = 
\deg_G^c(y) - \sigma_{z, X^+}(y)  
\stackrel{\text{Prop.\ref{prop:2}}}{\geq}    
\deg_G^c(y) - 3\ell 
\geq \delta^c(G) - 3\ell.  
\end{equation}
Then 
$$
\Delta(G) \geq \deg_G(y) = 
\sum_{\alpha \in \cA} |N_{\alpha}(y)|
+ 
\sum_{\beta \in \cB} |N_{\beta}(y)|
\stackrel{\eqref{eqn:7.19.2019.3:19p}}{\geq}    
\deg_F^-(y) 
+ 
\sum_{\beta \in \cB} |N_{\beta}(y)|  
\stackrel{\eqref{eqn:7.18.2019.6:46p}}{\geq}    
\deg_F^-(y) 
+ 
\delta^c(G) - 3\ell.   
$$
Corollary~\ref{cor:6} concludes the proof: since $\delta^c(G) \geq (n+1)/2$ holds by hypothesis, $\Delta(G)$ must satisfy
$$
\deg_F^-(y) \leq \Delta(G) - \delta^c(G) + 3\ell 
\stackrel{\text{Cor.\ref{cor:6}}}{<}  
\delta^c(G) + 4R + 3\ell - \delta^c(G) + 3\ell
\leq 
4R + 6\ell.   
$$

\end{document}